\documentclass[11pt,reqno]{amsart}

\usepackage{amsmath, amsfonts, amsthm, amssymb, color, cite}

\newtheorem{theorem}{Theorem}[section]
\newtheorem{lemma}{Lemma}[section]

\theoremstyle{definition}
\newtheorem{Definition}{Definition}[section]

\theoremstyle{remark}
\newtheorem{remark}{Remark}[section]

\numberwithin{equation}{section}
\allowdisplaybreaks

\newcommand{\R}{{\mathbb R}}

\def\f{\frac}

\def\hf1{^\f{1}{1-\xi^2}}

\def\be{\begin{equation}}
\def\en{\end{equation}}
\def\bs{\begin{split}}
\def\es{\end{split}}
\def\ba{\begin{align}}
\def\ea{\end{align}}

\newcommand{\eps}{\varepsilon}

\author[Feimin Huang]{Feimin Huang}
\address{School of Mathematical Sciences, University of Chinese Academy of Sciences, Beijing 100049, China; Academy of Mathematics and System Sciences, Chinese Academy of Sciences, Beijing 100190, China.}
\email{fhuang@amt.ac.cn}

\author[Tianhong Li]{Tianhong Li}
\address{Hua Loo-Keng Key Laboratory of Mathematics,
AMSS, CAS, Beijing 100190, China.}
\email{thli@math.ac.cn}

\author[Difan Yuan]{Difan Yuan}
\address{School of Mathematical Sciences, University of Chinese Academy of Sciences, Beijing 100049, China; Institute of Applied Mathematics,
 AMSS, CAS, Beijing 100190, China.}
\email{yuandf@amss.ac.cn}

\title[spherical symmetry]
{Global Entropy Solutions to Multi-Dimensional Isentropic Gas Dynamics with Spherical Symmetry
}


\date{\today}

\begin{document}
\begin{abstract}
We are concerned with spherically symmetric solutions to the Euler equations for the multi-dimensional compressible fluids, which have many applications in diverse real physical situations. The system can be reduced to one dimensional isentropic gas dynamics with geometric source terms. Due to the presence of the singularity at the origin, there are few papers devoted to this problem. The present paper proves two existence theorems of global entropy solutions. The first one focuses on the case excluding the origin in which the negative velocity is allowed, and the second one is corresponding to the case  including the origin with non-negative velocity. The $L^\infty$ compensated compactness framework and vanishing viscosity method are applied to prove the convergence of approximate solutions. In the second case, we show that if the blast wave initially moves outwards and the initial densities and velocities decay to zero with certain rates near origin, then the densities and velocities tend to zero with the same rates near the origin for any positive time. In particular, the entropy solutions in two existence theorems are uniformly bounded with respect to time.
\end{abstract}

\maketitle
\section{Introduction}

In this paper, we consider the Euler equations for compressible isentropic fluids with spherical symmetry which read
\begin{eqnarray}\label{sphericalsym}
\left\{ \begin{array}{ll}
\displaystyle \rho_t+\nabla\cdot \vec{m}=0,\,\vec{x}\in\R^{N},\\
\displaystyle \vec{m}_t+\nabla\cdot\left(\frac{\vec{m}\otimes\vec{m}}{\rho}\right)+\nabla p=0,\vec{x}\in\R^{N},
\end{array}
\right.
\end{eqnarray}
where $\rho, m$ and $p(\rho)$ denote the density, momentum and  pressure of the gas respectively. The pressure takes the form of $p(\rho)=p_0\rho^\gamma,$ with $p_0=\frac{\theta^{2}}{\gamma},\theta=\frac{\gamma-1}{2}$ and $\gamma>1$ being the adiabatic exponent.

We are interested in spherically symmetric solutions to system \eqref{sphericalsym} with the form
\begin{equation}\label{transform}
(\rho,\vec{m})(\vec{x},t)=(\rho(x,t),m(x,t)\frac{\vec{x}}{x}), x=|\vec{x}|.
\end{equation}
Then $(\rho(x,t),m(x,t))$ in \eqref{transform} is governed by the one-dimensional Euler equations with geometric source terms:
\begin{eqnarray}\label{geometric}
\left\{ \begin{array}{ll}
\displaystyle \rho_t+m_x=-\frac{N-1}{x}m,\,x\ge 0,t>0,\\
\displaystyle m_t+\left(\frac{m^2}{\rho}+p(\rho)\right)_x=-\frac{N-1}{x}\frac{m^{2}}{\rho},
x\ge 0,t>0.
\end{array}
\right.
\end{eqnarray}
For the system \eqref{geometric}, the study of spherically symmetric motion originates from several important applications such as the theory of explosion waves in medium, and the stellar dynamics including gaseous star formation and supernova formation. Note that the geometric source terms of \eqref{geometric} are singular at the origin, i.e., $x=0$. In the present paper, we first study the system \eqref{geometric} for the case that the origin is excluded. For simplicity, we consider \eqref{geometric} in the region outside the unit ball, that is,
\begin{eqnarray}\label{geometricc}
\left\{ \begin{array}{ll}
\displaystyle&v_t+F(v)_x=G(x,v),x\in(1,+\infty),t\in[0,+\infty),\\
\displaystyle &v|_{t=0}=v_0(x), x\in[1,+\infty),\\
&m|_{x=1}=0, t\in[0,+\infty),
\end{array}
\right.
\end{eqnarray}
with initial data $v_0(x)\in L^{\infty}([1,+\infty)),$ $v=(\rho, m)^\top, F(v)=(m, \frac{m^2}{\rho}+p(\rho))^\top,G(x,v)=$\\$
(a(x)m,a(x)\frac{m^{2}}{\rho})^\top,$ where $a(x)=-\frac{N-1}{x}.$
Then we consider  the case that the origin is included, i.e., $x\ge 0$. Note that in this case, the initial boundary value problem is equivalent to the Cauchy problem of compressible Euler equations \eqref{sphericalsym} with spherically initial data
\begin{equation}\label{sphericalinitial2}
\begin{aligned}
&(\rho, m)|_{t=0}=(\rho_0(x), m_0(x))\in L^{\infty}([0,+\infty)),x\geq0.\\
\end{aligned}
\end{equation} The boundary condition $m=0$ is based on the following: for classical solutions without vacuum to \eqref{geometric},  $\frac{dx(t)}{dt}=u(x(t), t)$ defines a particle path $x(t)$, where $u=\frac{m}{\rho}$ is the velocity. Any two particle paths $x_1(t)$ and $x_2(t)$ preserves the mass within $[x_1(0), x_2(0)]$. Therefore the particle path starting from the boundary should stay on the boundary, that implies $u=0$, i.e., $m=0$ on the boundary. It does not matter whether the boundary is $x=1$ or $x=0$.

There has been considerable progress on the existence of global entropy solutions for one dimension. T. Nishida \cite{[nishida]} first proved the existence of large BV solutions for isothermal gas (i.e., $\gamma=1$) by Glimm Scheme. Nishida and Smoller \cite{[nishida-smoller]} further studied the isentropic case (i.e., $\gamma>1$)  under some restrictions on the initial data. Note that both works mentioned above consider the case of excluding vacuum.

If the initial values contain vacuum, Diperna \cite{Diperna} first proved the global existence of $L^\infty$ entropy solutions with large initial data by the theory of compensated compactness for $\gamma=1+\frac{2}{2n+1},$ where $n\ge 2 $ is any  integer. Subsequently, Ding et al.\,\cite{Ding} and Chen et al.\,\cite{chen} successfully extended the result to $\gamma\in(1, \frac{5}{3}]$. Lions et al. \cite{Lions2} and \cite{Lions1} treated the case $\gamma>\frac{5}{3}.$ Finally  Huang et al. \cite{HuangWang} solved the existence problem of $L^\infty$ entropy solutions for $\gamma=1$ through compensated compactness and  analytic extension method.

For the inhomogeneous case, that is, the right hand side of \eqref{sphericalsym} is not zero, Ding et al. \cite{Ding1} established a general framework to investigate the global existence of $L^\infty$ entropy solutions through the fractional step Lax-Friedrichs scheme. It should be noted that in the framework of \cite{Ding1}, the approximate solutions are only required to be uniformly bounded in the space $x$, but not in the time $t$. The $L^\infty$ norm of approximate solutions may increase with respect to time $t$. Later on, there have been extensive works on the inhomogeneous case, see \cite{Chen2,Chen1998, Marcati, Marcati2} and the references therein. To study the large time behavior of entropy solution, it is important to show the uniform bound of solutions independent of time $t$.

As shown in \eqref{geometric}, the multi-dimensional compressible Euler system with spherical symmetry  can be reduced to one dimensional isentropic gas dynamics with geometric source terms, which may have singularity at $x=0.$
One of the main features is the resonance interaction among the characteristic mode and the geometrical source one.
The local existence of spherically symmetric  solutions outside a solid ball centered at the origin was discussed by Makino et al.\cite{Makino}
by the fractional step Lax-Friedrichs scheme.
The global existence was first studied by Chen and Glimm \cite{Chen2}, then by Tsuge   for global existence with uniform estimates\cite{Tsuge2006}.
For the domain including the origin, Chen et al. \cite{Chen1997} proved a global existence theorem with large $L^\infty$ data having only non-negative initial velocity.
More interesting works can be found in Chen\cite{Chen1997}, Chen and Li\cite{Chen2003}, Li and Wang\cite{Li}, Tsuge \cite{Tsuge1}\cite{Tsuge2006}, Yang\cite{Yang1995}\cite{Yang1996} and references therein. See also Wang and Wang \cite{Wangzejun}\cite{Wangzejun2}.
 For more background of physical motivation for studying spherically symmetric solutions, please refer to \cite{Courant,Dafermos}. Recently Chen et al.\cite{ChenP} established  $L^p$ global finite-energy entropy solution of the isentropic Euler equations with spherical symmetry and large initial data.

Note that all of above works are either based on numerical schemes, which need laborious estimates, or related with special solutions. In this paper, we apply the vanishing viscosity method together with the invariant region of parabolic system with nonlinear source terms
to obtain  a priori uniform estimates of viscosity solutions. The revised version of the theory of invariant region (Lemma 2.1) is quite powerful and easy to be used in dealing with  source terms. This approach is valid for both of the Cauchy problem and initial boundary value problem. In the first part of present paper, we consider an initial boundary value problem outside a unit ball in which the origin is excluded. In the second part of present paper, we study the Cauchy problem so that the origin is included. In both cases, we obtain the uniform bound of viscosity solutions independent of   time $t$, while the $L^\infty$ bound depends on time $t$ in almost all previous works. This marks an important step to investigate the large time behavior of entropy solutions. It is worth pointing out that a set of new delicate control functions is designed to obtain the uniform bound of approximate solutions. Moreover, a new approach is proposed in the proof of lower bound of density by using the appropriate decomposition of source terms of heat equation and corresponding solutions successively.
For the case excluding the origin, we allow negative velocity when gas initially moves inwards. The major difficulties arise from handling  the solid ball $x=1$ and far field $x=\infty$ at the same time. For the case including the origin, we also observe a new phenomena that if the blast wave initially moves outwards and the initial densities and velocities decay to zero with certain rates near origin, then the densities and velocities tend to zero with the same rates near the origin for any positive time. It is remarked that in this case, the cavity phenomena occurs in the origin, but the gases move radially outwards. The method in our paper can be applied to solve the Euler equations with source terms, especially with geometric effect, such as two-dimensional radial gas flow, in gas flow through general nozzle, etc. All of these results will be discussed in a forthcoming paper.
Before formulating the main results, we define the entropy solutions of initial boundary value problem and Cauchy problem respectively as follows.

\begin{Definition}\label{def2}
A measurable function $v(x, t)\in L^\infty ([1,+\infty)\times \R^+)$ is called a global entropy solution of the initial boundary value problem \eqref{geometricc} provided that
\begin{equation}\label{weakexclude}
\int_0^{+\infty}\int_1^{+\infty}\left(v\Phi_t+F(v)\Phi_x+G(x,v)\Phi\right)dxdt+\int_1^{+\infty}v_0(x)\Phi(x, 0)dx=0
\end{equation}
holds for any function $\Phi\in C^1_0((1,+\infty)\times[0,+\infty))$ and
for any weak convex entropy pair $(\eta, q)(v)$, the inequality
\begin{equation}\label{entropyinequa}
\eta(v)_t+q(v)_x-\nabla\eta(v)\cdot G(x,v)\leq0
\end{equation}
holds in the sense of distributions.
\end{Definition}

\

The weak convex entropy-flux pair will be defined in the next section. The precise statement of the first result is given below.

\

\begin{theorem}\label{theorem1.1}\text(Excluding the origin)\label{exclude}
Let $1<\gamma\leq3.$ Given any positive constant $M_2$, there exists a constant $M_1$, which is bigger than $M_2$, such that if the initial and boundary data satisfy
 \begin{equation}\label{ini}
\begin{aligned}
&\rho_0(x)\ge 0,~\frac{m_0}{\rho_0}+\rho_0^{\theta} \leq M_1-M_2x^{-\alpha},\,\frac{m_0}{\rho_0}-\rho_0^{\theta}\geq -M_2x^{-\alpha},\text{ a.e.},\\
&\qquad\qquad\quad\quad m(1,t)=0,t>0,\\
\end{aligned}
\end{equation}
then there exists a global entropy solution of \eqref{geometricc}  satisfying
\begin{equation}\label{solution}
\begin{aligned}
\rho(x,t)\ge 0,~(\frac{m}{\rho}+\rho^{\theta})(x,t) \leq M_1-M_2x^{-\alpha},
\,(\frac{m}{\rho}-\rho^{\theta})(x,t)\geq -M_2x^{-\alpha},\text{ a.e.},
\end{aligned}
\end{equation} in the sense of Definition \ref{def2},
where  $\alpha$ is any constant satisfying $\frac{(N-1)\theta}{(1+\sqrt{\theta})^2}\leq\alpha \leq\frac
{(N-1)\theta}{(1-\sqrt{\theta})^2}.$
\end{theorem}
\begin{remark}For any $M>M_1$, Theorem 1.1 still holds if $M_1$ is replaced by $M$.
\end{remark}
\begin{remark}From (\ref{solution}), it implies that $0\le\rho^\theta\le \frac{M_1}{2}, -M_2x^{-\alpha}\le\frac{m}{\rho}\le M_1-M_2x^{-\alpha}$.
\end{remark}
\begin{remark}
From \eqref{ini} and \eqref{solution}, the negative velocity is allowed.
	\end{remark}
\begin{remark} The lower bound of $\alpha $ is much more important than its upper bound. When $\alpha$ becomes smaller, the increasing rate of $-x^{-\alpha}$ to $0$ decreases, which means that the range of negative velocity is larger.
	In Tsuge \cite{Tsuge2004a}, the initial data satisfies
	$\frac{m_0}{\rho_0}-\rho_0^{\theta}\geq -Cx^{-\frac{(N-1)\theta}{1+\theta}}$.
	Since $\frac{(N-1)\theta}{(1+\sqrt{\theta})^2}< \frac{(N-1)\theta}{1+\theta}$, more initial values with negative velocities are allowed here.
\end{remark}

\begin{Definition}\label{def3}
A measurable function $v(x, t)$ is called a global entropy solution of the Cauchy problem \eqref{geometric} and \eqref{sphericalinitial2} provided that
\begin{equation}\label{weakinclude}
\int_0^{+\infty}\int_0^{+\infty}\left(v\Phi_t+F(v)\Phi_x+G(x,t,v)\Phi\right)dxdt+\int_0^{+\infty}v_0(x)\Phi(x, 0)dx=0
\end{equation}
holds for any function $\Phi\in C^1_0((0,+\infty)\times[0,+\infty))$ and
for any weak convex entropy pair $(\eta, q)$, the inequality
\begin{equation*}
\eta(v)_t+q(v)_x-\nabla\eta(v)\cdot G(x,t,v)\leq0
\end{equation*}
holds in the sense of distributions.
\end{Definition}
The second result of this paper is stated as follows.
\begin{theorem}\text(Including the origin)\label{include}
Let $\gamma>1.$ Assume that for any nonnegative constants $c$ and $M_3,$   there hold
 \begin{equation}\label{ini2}
\begin{aligned}
&\rho_0(x)\ge 0,~\frac{m_0}{\rho_0}+\rho_0^{\theta} \leq M_3x^{c\theta},\,\frac{m_0}{\rho_0}-\rho_0^{\theta}\geq 0,\text{ a.e. $x\in [0,+\infty)$},\\
\end{aligned}
\end{equation}
then there exists a global entropy solution of \eqref{geometric} and \eqref{sphericalinitial2} satisfying
\begin{equation}\label{solution2}
0\leq\rho(x,t)\leq (\frac{M_3}{2})^{\frac{1}{\theta}}x^c,~0\leq m(x, t)\leq M_3\rho(x, t)x^{c\theta} \text{ a.e.$(x,t)\in[0,+\infty)\times\R^{+}$}.
\end{equation}
\end{theorem}

\begin{remark}
	Theorem \ref{include} means that if the blast wave initially moves outwards and $\rho^\theta$ and  $u=\frac{m}{\rho}$ initially decay to zero with certain rates near origin, then they tend to zero with same rates near the origin for any positive time.
	\end{remark}
\begin{remark}
	In Theorem \ref{include}, the initial data can be allowed to tend to infinity  at far field.
\end{remark}
\begin{remark}
The invariant region $\frac{m}{\rho}+\rho^\theta\le M, \frac{m}{\rho}-\rho^\theta\ge 0$ was first observed in \cite{Chen1997}.  This  corresponds to the speical case  $c=0$ in Theorem \ref{include}.
\end{remark}

The main ingredient in proving Theorem \ref{theorem1.1} is how to get the uniform estimates of viscosity solutions independent of  viscosity $\varepsilon$ and time $t$. In fact, the viscosity solutions are uniformly bounded through a maximum principle for parabolic system, see Lemma \ref{initial maximum} below. Roughly speaking, we first add viscous perturbation to the system \eqref{geometric} and get a viscous system \eqref{geometricvis}, which can be reduced into a decoupled system \eqref{wwzz} of the Riemann invariants.
Unfortunately Lemma \ref{initial maximum} can not be directly applied since the coefficients $a_{12}$ and $a_{21}$ of Lemma \ref{initial maximum} may not be negative in the system \eqref{wwzz}, while they have to be negative for the application of the  lemma.  The key point is to introduce modified Riemann invariants to derive a new system \eqref{rst1} in which the coefficients have desired sign. To prove Theorem \ref{include}, we first introduce a space scaling transformation for both variables $\rho=\tilde{\rho}x^c,m=\tilde{m}x^d$ and space coordinate $\xi=\frac{1}{c-d+1}x^{c-d+1}\quad(\mbox{if}~c-d+1\ne 0)$  or $\xi=\ln x \quad(\mbox{if}~c-d+1=0)$, see Section 4 below.  Then we add viscous perturbations $(\varepsilon \tilde{\rho}_{\xi\xi}, \varepsilon \tilde{m}_{\xi\xi})$ to the new system for $(\tilde{\rho},\tilde{m})$. Again using Lemma \ref{initial maximum}, we get desired uniform estimates of viscosity solutions being independent of viscosity $\varepsilon$ and time $t$.


\

The present paper is organized as follows: In Section \ref{formula}, we  construct approximate solutions by adding viscosity and some preliminaries are given. In Section \ref{excludeorigin}, we first obtain the uniform upper bounds independent of viscosity $\varepsilon$ and time $t$ for the viscosity solutions and then prove the $H_{loc}^{-1}$ compactness for  entropy-entropy flux pairs, and finally Theorem \ref{theorem1.1}. Section \ref{includeorigin} is devoted to the proof of Theorem \ref{include}, i.e., the existence of entropy solutions with spherical symmetry.

\section{Preliminaries and Formulations}\label{formula}
We first introduce some basic facts for the system \eqref{geometric}.
The eigenvalues are
\begin{equation}\label{2.1}
\lambda_1=\frac{m}{\rho}-\theta\rho^\theta,\quad
\lambda_2=\frac{m}{\rho}+\theta\rho^\theta,\quad
\end{equation}
where $\theta=\frac{\gamma-1}{2},$\,and the corresponding right eigenvectors are
\begin{equation}\label{2.2}
r_1=\left[\begin{array}{cc}
1\\ \lambda_1
\end{array}
\right],\quad
r_2=\left[\begin{array}{cc}
1\\ \lambda_2
\end{array}
\right].
\end{equation}
The Riemann invariants $(w, z)$ are given by
\begin{equation}\label{2.3}
w=\frac{m}{\rho}+\rho^\theta,\quad z=\frac{m}{\rho}-\rho^\theta,
\end{equation}
satisfying $\nabla w\cdot r_1=0$ and $\nabla z\cdot r_2=0,$ where $\nabla=(\partial_\rho, \partial_m)$ is the gradient with respect to $U=(\rho,m)$. A pair of functions $(\eta, q): \R^+\times\R\mapsto\R^2$ is defined to be an entropy-entropy flux pair of system  \eqref{geometric} if it satisfies
\begin{equation}\label{2.4}
\nabla q(U)=\nabla\eta(U)\nabla\left[\begin{array}{ccc}
m\\ \frac{m^2}{\rho}+p(\rho)
\end{array}
\right].
\end{equation}
When $$\eta\left|_{ \frac{m}{\rho}\text{ fixed }}\rightarrow 0,\qquad\mbox{as }\rho\rightarrow 0, \right.$$ $\eta(\rho, m)$ is called weak entropy.
Moveover, an entropy $\eta(\rho,m)$ is  convex (strictly convex) if the Hessian matrix $\nabla^2\eta(\rho, m)$ is nonnegative (positive). For example,
\begin{equation}\label{2.5}
\eta^*(\rho, m)=\frac{m^2}{2\rho}+\frac{p_0\rho^\gamma}{\gamma-1},~~
q^*(\rho, m)=\frac{m^3}{2\rho^2}+\frac{\gamma p_0\rho^{\gamma-1}m}{\gamma-1},
\end{equation}
is a strictly convex entropy pair. 
As shown in \cite{Lions2} and \cite{Lions1},  any weak entropy for the system \eqref{geometric}  is
\begin{equation}\label{2.6}
\begin{split}
\eta=\rho\int_{-1}^1g(\frac{m}{\rho}+\rho^\theta s)(1-s^2)^\lambda ds,~~
q=\rho\int_{-1}^1(\frac{m}{\rho}+\rho^\theta\theta s)g(\frac{m}{\rho}+\rho^\theta s)(1-s^2)^\lambda ds,
\end{split}
\end{equation}
with  $\lambda=\frac{3-\gamma}{2(\gamma-1)}$ and $g(\cdot)\in C^2(\R)$ is any function.\\

Now, we will introduce a revised version of the theory of invariant region which is essentially based on the maximum principle for parabolic equation, see \cite{Smoller}. 
\begin{lemma}(Maximum principle on bounded domain)\label{initial maximum}
Let $p(x,t), q(x,t)$, $(x,t)\in[a,b]\times[0,T]$ be any bounded classical solutions of the following quasilinear parabolic system
\begin{eqnarray}\label{pq}
\left\{ \begin{split}
\displaystyle &p_t+\mu_1 p_x=
p_{xx}+a_{11}p+a_{12}q+R_1,\\
\displaystyle &q_t+\mu_2 q_x=
q_{xx}+a_{21}p+a_{22}q+R_2,
\end{split}
\right.
\end{eqnarray}
with
\[\begin{aligned}
  p(x,0)&\leq 0,~q(x, 0)\geq0, \text{ for } ~x\in[a,b],\\
  p(a,t)& \leq 0,\,q(a,t)\geq0, \text{ for }~ t\in[0,T],\\
  p(b,t)& \leq 0,\,q(b,t)\geq0, \text{ for }~ t\in[0,T],\\
\end{aligned}\]
where $$\mu_{i}=\mu_i(x,t,p(x,t),q(x,t)),a_{ij}=a_{ij}(x,t,p(x,t),q(x,t))$$ and the source terms $$R_i=R_i(x,t,p(x,t),q(x,t),p_{x}(x,t),q_{x}(x,t)),i,j=1,2,\forall(x,t)\in[a,b]\times[0,T].$$ Here $\mu_{i},a_{ij}$ are bounded with respect to $(x,t,p,q)\in[a,b]\times[0,T]\times K,$ where $K$ is an arbitrary compact subset in $\R^2.$  $a_{12},a_{21},R_{1},R_{2}$ are continuously differentiable with respect to $p,q.$ Assume that
\begin{description}
	\item[(C1)]  $a_{12}\leq0$ holds for $p=0$ and $q\geq0$;  $a_{21}\leq0$ holds for $q=0$ and $p\leq0$;\label{jj}
	\item[(C2)]  $~R_1\leq0$ holds for $p=0$ and $q\geq0$; $R_2\geq0$ holds for $q=0$ and $p\leq0$.
\end{description}
Then for any $(x, t)\in[a,b]\times[0,T],$
$$p(x,t)\leq 0, ~~q(x, t)\geq0.$$
\end{lemma}
\section{Proof of Theorem 1.1}\label{excludeorigin}

\subsection{Uniform upper bound estimate}\label{upper}
We approximate \eqref{geometric} by adding artificial viscosity as follows:
\begin{eqnarray}\label{geometricvis}
\left\{ \begin{array}{ll}
\displaystyle \rho_t+m_x=-\frac{N-1}{x}m+\eps \rho_{xx},\,\\
\displaystyle m_t+\left(\frac{m^2}{\rho}+p(\rho)\right)_x=-\frac{N-1}{x}\frac{m^{2}}{\rho}+\eps m_{xx}-\frac{2\eps \alpha M_2}{x^{\alpha+1}}\rho_x.
\end{array}
\right.
\end{eqnarray}
We consider \eqref{geometricvis} on a cylinder $(1,b)\times \R^{+}$, with $\R^{+}=[0,+\infty),  b:=b(\eps)$ satisfying $\lim\limits_{\eps\rightarrow0}b(\eps)=\infty.$ The initial-boundary values are given by
\begin{equation}\label{geometricini-vis}
\begin{aligned}
&(\rho, m)|_{t=0}=(\rho_0^\eps(x), m_0^\eps(x))=(\rho_0(x)+\eps^{\frac{2}{\theta}}, \frac{m_0(x)}{\rho_0(x)}(\rho_0(x)+\eps^{\frac{2}{\theta}}))\ast j^\eps, 1\leq x\leq b,\\
&~~~~~~~~~~(\rho,m)|_{x=1}=(\rho_0^\eps(1),0),(\rho,m)|_{x=b}=(\rho_0^\eps(b),m_0^\eps(b)),  t>0,
\end{aligned}
\end{equation}
where $j^\eps$ is a standard mollifier with small parameter $\eps>0$. 
Next we will derive the uniform bound of the viscosity solutions by the maximum principle, i.e., Lemma \ref{initial maximum}.
Based on the Riemann invariants $w$ and $z$, we transform \eqref{geometricvis} into the following form:
\begin{eqnarray}\label{wwzz}
\left\{ \begin{split}
\displaystyle w_t+\lambda_2 w_x=&\eps w_{xx}+2\eps(w_x-\alpha M_2 x^{-\alpha-1})\frac{\rho_x}{\rho}-\theta \frac{N-1}{x}\rho^{\theta}z-\theta \frac{N-1}{x}\rho^{2\theta}\\
&-\eps\theta(\theta+1)\rho^{\theta-2}\rho_x^2,\\
\displaystyle z_t+\lambda_1 z_x=&\eps z_{xx}+2\eps (z_x-\alpha M_2 x^{-\alpha-1})\frac{\rho_x}{\rho}+\theta \frac{N-1}{x}\rho^{\theta}z+\theta \frac{N-1}{x}\rho^{2\theta}\\
&+\eps\theta(\theta+1)\rho^{\theta-2}\rho_x^2.
\end{split}
\right.
\end{eqnarray}
From the initial condition \eqref{ini}, we set the  control functions $(\phi,\psi)$:
$$\phi=M_1-M_2x^{-\alpha}+\eps e^{Ct},\,\psi=M_2x^{-\alpha}+\eps e^{Ct}, $$ 
where $C$ and $M_1$  are positive constants which will be determined later.
Then a simple calculation shows that
\begin{equation*}
\begin{split}
&\phi_t=\eps C e^{Ct},~\phi_x=\alpha M_2x^{-\alpha-1}, ~\phi_{xx}=-\alpha(\alpha+1)M_2x^{-\alpha-2};\\
&\psi_t=\eps C e^{Ct},~\psi_x=-\alpha M_2 x^{-\alpha-1}, ~\psi_{xx}=\alpha(\alpha+1)M_2x^{-\alpha-2}.
\end{split}
\end{equation*}
Define a modified Riemann invariants $(\bar{w},\bar{z})$ as
\begin{equation}\label{r}
\bar{w}=w-\phi, ~~\bar{z}=z+\psi.
\end{equation}
We shall use Lemma \ref{initial maximum} to show $\bar{w}\le 0$ and $\bar{z}\ge 0$ for any time.
Inserting \eqref{r} into \eqref{wwzz} yields the equations for $\bar{w}$ and $\bar{z}:$
\begin{eqnarray}\label{phipsi1}
\left\{\begin{split}
\bar{w}_t+\lambda_2\bar{w}_x=&\eps\bar{w}_{xx}
+2\eps\frac{\rho_x}{\rho}\bar{w}_x-\eps\theta(\theta+1)\rho^{\theta-2}\rho_x^2\\
&-\phi_t-\lambda_2\phi_x+\eps\phi_{xx}-\theta \frac{N-1}{x}\rho^{\theta}(\bar{z}-\psi)-\theta \frac{N-1}{x}\rho^{2\theta}, \\
\bar{z}_t+\lambda_1\bar{z}_x=&\eps\bar{z}_{xx}
+2\eps\frac{\rho_x}{\rho}\bar{z}_x+\eps\theta(\theta+1)\rho^{\theta-2}\rho_x^2\\
&+\psi_t+\lambda_1\psi_x-\eps\psi_{xx}+\theta \frac{N-1}{x}\rho^{\theta}(\bar{z}-\psi)+\theta \frac{N-1}{x}\rho^{2\theta}.\\
\end{split}
\right.
\end{eqnarray}
Note that $$\lambda_1=\bar{z}-\psi+(1-\theta)\rho^\theta,\lambda_2=\bar{w}+\phi+(\theta-1)\rho^\theta.$$
The system \eqref{phipsi1} becomes
\begin{eqnarray}\label{rst1}
\displaystyle\left\{ \begin{split} &\bar{w}_t+(\lambda_2-2\eps\frac{\rho_x}{\rho})\bar{w}_x
=\eps\bar{w}_{xx}+a_{11}\bar{w}
+a_{12}\bar{z}+R_1,\\
&\bar{z}_t+(\lambda_1-2\eps\frac{\rho_x}{\rho})\bar{z}_x
=\eps\bar{z}_{xx}+a_{21}\bar{w}
+a_{22}\bar{z}+R_2,
\end{split}
\right.
\end{eqnarray}
where
\begin{equation*}
\begin{split}
&a_{11}=-\phi_x,\,a_{12}=-\theta \frac{N-1}{x}\rho^{\theta}\leq0,\,a_{21}=0,\,a_{22}=\psi_x+\theta\frac{N-1}{x} \rho^{\theta},\\
\end{split}
\end{equation*}
$$R_1=-\phi_t-[\phi+(\theta-1)\rho^\theta]\phi_x+\eps\phi_{xx}-\eps\theta(\theta+1)\rho^{\theta-2}\rho_x^2+\theta\frac{N-1}{x}\rho^{\theta}\psi-\theta\frac{N-1}{x}\rho^{2\theta},$$
$$R_2=\psi_t+[-\psi+(1-\theta)\rho^\theta]\psi_x-\eps\psi_{xx}+\eps\theta(\theta+1)\rho^{\theta-2}\rho_x^2-\theta\frac{N-1}{x}\rho^{\theta}\psi+\theta\frac{N-1}{x}\rho^{2\theta}.$$
A direct computation gives
\begin{equation*}
\begin{split}
R_1\leq&(M_1-M_2x^{-\alpha}+2\eps e^{Ct}+(\theta-1)\rho^\theta)\psi_x+\theta\frac{N-1}{x}\rho^\theta\psi-\theta\frac{N-1}{x}\rho^{2\theta}-\eps Ce^{Ct}\\
=&-\alpha M_1M_2x^{-\alpha-1}+\alpha M_2^2x^{-2\alpha-1}+\alpha M_2(1-\theta)\rho^\theta x^{-\alpha-1}+\theta(N-1)M_2\rho^\theta x^{-\alpha-1}\\
-&\theta\frac{N-1}{x}\rho^{2\theta}-2\eps e^{Ct}M_2x^{-\alpha-1}-\eps Ce^{Ct}+\theta\frac{N-1}{x}\rho^{\theta}\eps e^{Ct}.
\end{split}
\end{equation*}
To ensure $R_1\leq0,$ it is sufficient to show that
\begin{equation}\label{37}
\begin{split}
&\alpha M_2^2x^{-\alpha}+M_2\left[\alpha(1-\theta)+\theta(N-1)\right]\rho^\theta+
\theta(N-1)\rho^{\theta}\eps e^{Ct}x^\alpha\\ \leq&\alpha M_1M_2+
\theta(N-1)\rho^{2\theta} x^\alpha+2\eps e^{Ct}M_2+\eps Ce^{Ct}x^{\alpha+1}.
\end{split}
\end{equation}
For the second term on the left hand side of \eqref{37}, by Cauchy-Schwartz inequality, we have
$$M_2\left[\alpha(1-\theta)+\theta(N-1)\right]\rho^\theta\leq\frac{1}{2}\theta(N-1)
\rho^{2\theta}+\frac{M_2^2\left[\alpha(1-\theta)+\theta(N-1)\right]^2}{2\theta(N-1)}.$$
Since $x\ge 1$, the first and second terms in \eqref{37} can be controlled by choosing $M_1\ge M_2$ and large enough.

For the third term,  since $\rho^{\theta}\leq\frac{1}{2}(1+\rho^{2\theta})$,  we choose $C$ being at least bigger than $\frac{\theta(N-1)}{2}$ and we will see that  the choice of $C$ also depends on $R_2$ later. Then $\theta(N-1)\rho^{\theta}\eps e^{Ct}x^\alpha$ can be controlled by $\eps Ce^{Ct}x^{\alpha+1}$ and $\theta(N-1)\rho^{2\theta} x^\alpha$,
as long as for fixed $T$, choosing $\epsilon$ small enough such that $\sqrt{\eps} e^{Ct}\le \sqrt{\eps} e^{CT}<1$ for $t\le T$. So it can be seen that  $\epsilon$ relies on $C$ and $T$.

 Thus $R_1\le 0$.

We now turn to the term $R_2$. Denote $\beta=\frac{\alpha}{\theta(N-1)}$. By a direct calculation, we have
\begin{equation}\label{38}
\begin{split}
R_2\geq&\left[-\psi+(1-\theta)\rho^\theta\right]\psi_x-\eps\psi_{xx}-\theta\frac{N-1}{x}\rho^\theta\psi+\theta\frac{N-1}{x}\rho^{2\theta}+\psi_t\\
=&\frac{\theta(N-1)}{x}\left[\rho^{2\theta}-\frac{\alpha(1-\theta)+\theta(N-1)}{\theta(N-1)}\rho^\theta\psi+\frac{\alpha\psi^2}{\theta(N-1)}\right]\\
-&\eps \alpha\psi e^{Ct}\frac{1}{x}+\eps(1-\theta)\rho^{\theta}\alpha e^{Ct}\frac{1}{x}+\eps Ce^{Ct}-\eps \alpha(\alpha+1)M_2x^{-\alpha-2}\\
=&\frac{\theta(N-1)}{x}\left[\left(\rho^\theta-\frac{\beta(1-\theta)+1}{2}\psi\right)^2-\frac{(\beta(1-\theta)+1)^2}{4}\psi^2+\beta\psi^2\right]\\
+&\eps\left[- \alpha(M_2x^{-\alpha}+\eps e^{Ct}) e^{Ct}\frac{1}{x}+ Ce^{Ct}- \alpha(\alpha+1)M_2x^{-\alpha-2}+(1-\theta)\rho^{\theta}\alpha e^{Ct}\frac{1}{x}\right].\\
\end{split}
\end{equation}
For the second term of right hand side of \eqref{38}, choosing $C$ also sufficiently large,  and  using $\sqrt{\eps} e^{Ct} <1$ from the estimate of $R_1$ above, then its positivity is obtained.
To ensure $R_2\geq0,$ it remains to guarantee that $$\beta\psi^2-\frac{(\beta(1-\theta)+1)^2}{4}\psi^2\geq0,$$
i.e.,$$g(\beta):=\beta^2(1-\theta)^2-2\beta(1+\theta)+1\leq0,$$
which holds for any
 $$\theta(N-1)\beta_1\leq\alpha\leq\theta(N-1)\beta_2, $$
 where $$\beta_1=\frac{1}{(1+\sqrt{\theta})^2}, ~\beta_2=\frac{1}{(1-\sqrt{\theta})^2},$$
 are the roots of the equation $g(\beta)=0.$
By the initial and boundary data \eqref{geometricini-vis}, we obtain $$\bar{w}(x,0)=w(x,0)-\phi\leq 0,\,\bar{z}(x, 0)=z(x, 0)+\psi\geq0,\,\text{for }\, x\geq 1;$$
$$\bar{w}(1,t)\leq 0,\,\bar{z}(1,t)\geq0,\,\text{for } t>0;$$
$$\bar{w}(b,t)\leq 0,\,\bar{z}(b,t)\geq0,\,\text{for } t>0.$$
By Lemma \ref{initial maximum}, $$w(x,t)\leq M_1-M_2x^{-\alpha}+\eps e^{Ct}, z(x,t)\geq -M_2x^{-\alpha}-\eps e^{Ct},\qquad\mbox{i.e.,}$$
$$w(x,t)\leq M_1-M_2x^{-\alpha}+\sqrt{\eps}, z(x,t)\geq -M_2x^{-\alpha}-\sqrt{\eps},$$
where we have used the fact that  $\sqrt{\eps} e^{Ct}<1.$
 This gives the following Theorem \ref{main3.1} for approximate solutions.
\begin{theorem}\label{theorem 3.1}\text($L^{\infty}$ estimate: excluding the origin)\label{main3.1}  Let $1<\gamma\leq3.$ Given any positive constant $M_2$, there exists a constant $M_1$, which is bigger than $M_2$,   such that if
\begin{equation}\label{ini1}
\begin{aligned}
&\rho_0(x)\ge \eps^{\frac{2}{\theta}} ,~\frac{m_0}{\rho_0}+\rho_0^{\theta} \leq M_1-M_2x^{-\alpha}+\eps,\, \frac{m_0}{\rho_0}-\rho_0^{\theta}\geq -M_2x^{-\alpha}-\eps,\text{$x\in [1,b]$},\\
&\left(\frac{m}{\rho}+\rho^\theta\right)|_{x=1}\leq M_1-M_2 ,~\left(\frac{m}{\rho}-\rho^\theta\right)|_{x=1}\geq -M_2,\\
&\left(\frac{m}{\rho}+\rho^\theta\right)|_{x=b}\leq M_1-M_2b^{-\alpha}+\eps,~\left(\frac{m}{\rho}-\rho^\theta\right)|_{x=b}\geq-M_2b^{-\alpha}-\eps,\,t>0,
\end{aligned}
\end{equation}
where $\frac{(N-1)\theta}{(1+\sqrt{\theta})^2}\leq\alpha\leq\frac{(N-1)\theta}{(1-\sqrt{\theta})^2}$, then the solutions of \eqref{geometricvis}-\eqref{geometricini-vis} satisfy
\begin{equation}
(\frac{m}{\rho}+\rho^{\theta})(x,t) \leq M_1-M_2x^{-\alpha}+\sqrt{\eps},
\,(\frac{m}{\rho}-\rho^{\theta})(x,t)\geq -M_2x^{-\alpha}-\sqrt{\eps},
\end{equation}
\end{theorem}
The only thing left to be checked in Theorem \ref{main3.1} is the condition (\ref{ini1}) for initial data. It is easy to see that the function $f(r)=(\eps+r^\theta)^{\frac{1}{\theta}}-r-\eps^{\frac{2}{\theta}}$ is  increasing on  $(0, \infty)$ when $\theta\in (0, 1]$ and $\eps$ is small. By $f(0)>0$, it follows that $f(r)>0$  for $r\ge 0$. Hence the condition (\ref{ini1}) holds.
\
\subsection{Lower bound estimate of density}\label{low}
From the above argument, we know that the velocity $u=\frac{m}{\rho}$ is uniformly bounded, i.e., $|u|\le M_1$. The lower bound of density can be derived as in \cite{Lu}, but we introduce a different treatment for estimating a \textit{priori} lower bound for solutions of heat equations with general source terms. Set $e=\ln\rho$, and we can get a scalar equation for $e$, that is,
\begin{equation}\label{v}
e_t+e_xu+u_x=\eps e_{xx}+\eps e_x^2-\frac{N-1}{x}u.
\end{equation}
From the initial-boundary value of \eqref{geometricini-vis},  we have $e|_{t=0}=\ln\rho_0^\eps(x),e|_{x=1}=\ln\rho_0^\eps(1),e|_{x=b}=\ln\rho_0^\eps(b)$.
Then it follows from \eqref{v} that
$$e_t-\eps e_{xx}=\eps(e_x-\frac{u}{2\eps})^2-\frac{u^2}{4\eps}-\frac{N-1}{x}u-u_x.$$
To the best of our knowledge, we found no literature that ever stated the following Lemma to obtain the lower bound of density.
\begin{lemma}\label{density}
	Assume that $w$ is a classical solution of heat equation:
	\begin{eqnarray}\label{u}
	\left\{ \begin{array}{ll}
	\displaystyle w_t-\eps w_{xx}=f_1(x,t)+f_2(x,t)+h_x(x,t),a<x<b,t>0,\,\\
	\displaystyle w|_{t=0}=\varphi(x), a<x<b,\\
	\displaystyle w|_{x=a}=\varphi(a),w|_{x=b}=\varphi(b),
	\end{array}
	\right.
	\end{eqnarray}
	where $f_1(x,t)\geq0, f_2(x,t)$ and $h(x,t)$ are bounded smooth functions and $\varphi(x)$ is a bounded function. Then, there exists a positive constant $C(a,b,\eps,t)$, s.t., $$w\geq-C(a,b,\eps,t).$$
\end{lemma}
\noindent Proof: the proof is given in Appendix.

\

By applying Lemma \ref{density}, for $$w=e, f_1(x,t)=\eps(e_x-\frac{u}{2\eps})^2, f_2(x,t)=-\frac{u^2}{4\eps}-\frac{N-1}{x}u, h(x,t)=u,$$ it follows that $$\rho\geq e^{-C(a,b,\eps, t)}>0.$$ 

From the local existence of approximate  solutions, the upper and the lower bound of density, we can conclude the following theorem for the global existence of approximate solutions.
\begin{theorem}\label{theorem 3.2} Under the assumption of the previous theorem,
	for any time $T>0,$ there exists $\eps_0$  such that for $0<\eps<\eps_0$, the initial-boundary value problem
	\eqref{geometricvis}-\eqref{geometricini-vis} admits a unique classical solution on $[1, b]\times [0, T]$ satisfying
\begin{equation*}
\rho(x, t)\ge e^{-C(\eps, T)}, (\frac{m}{\rho}+\rho^{\theta})(x,t) \leq M_1-M_2x^{-\alpha}+\sqrt{\eps},
\,(\frac{m}{\rho}-\rho^{\theta})(x,t)\geq -M_2x^{-\alpha}-\sqrt{\eps}.
	\end{equation*}
\end{theorem}
\subsection{$H^{-1}_{loc}$ compactness of the entropy pair}\label{entropy11}
For any $T\in(0, \infty)$, let $\Pi_T=(1,+\infty)\times (0, T).$  We consider the entropy dissipation measures
\begin{equation}\label{h-1}
\eta(\rho, m)_t+q(\rho, m)_x,
\end{equation}
where $(\eta,q)$ is  any weak entropy-entropy flux pair whose formula is given in \eqref{2.6}. We will apply the Murat Lemma to conclude that the entropy dissipation measures in \eqref{h-1} lie in a compact set of $H^{-1}_{loc}(\Pi_T)$ .
 \begin{lemma}{(Murat \cite{Murat})}\label{murat}
	Let $\Omega\subseteq\R^n$ be an open set, then
	\begin{equation*}
		(\text{compact set of } W^{-1, q}_{loc}(\Omega))\cap(\text{bounded set of } W^{-1, r}_{loc}(\Omega))\\
		\subset(\text{compact set of } H^{-1}_{loc}(\Omega)),
	\end{equation*}
	where $1<q\leq 2<r.$
\end{lemma}
We write $U=-\frac{N-1}{x}m, V=-\frac{N-1}{x}\frac{m^2}{\rho}$ for simplicity. Let $K\subset\Pi_T$ be any compact set and choose
$\varphi\in C_c^\infty(\Pi_T)$ such that $\varphi|_{K}=1$ and $0\leq\varphi\leq1$. When $\eps$ is small, $K\subset (1,b(\eps))\times (0, T)$.
Multiplying \eqref{geometricvis} by  $\nabla\eta^* \varphi$ with $\eta^*$ being the mechanical entropy in \eqref{2.5},
we obtain
\begin{equation}\label{4.1}
\begin{split}
&\eps\int\int_{\Pi_T}(\rho_x, m_x)\nabla^2\eta^*(\rho_x, m_x)^\top \varphi dxdt\\
=&\int\int_{\Pi_T}\left[(V\phi-2\eps \phi_x\rho_x)\eta^*_m+U\eta^*_\rho\varphi
+\eta^*\varphi_t+q^*\varphi_x+\eps\eta^*\varphi_{xx}\right]dxdt.
\end{split}
\end{equation}
Note that
\[
(\rho_x, m_x)\nabla^2\eta^*(\rho_x, m_x)^\top=p_0\gamma\rho^{\gamma-2}\rho_x^2
+\rho u_x^2,
\]
and
\begin{equation*}
|(V\varphi-2\eps\phi_x\rho_x)\eta^*_m|\leq\frac{\eps p_0\gamma}{2}\rho^{\gamma-2}\rho_x^2+\eps C\phi_x^2m^2\rho^{-\gamma}+C\|V\|_{L^{\infty}},
\end{equation*}
we get
\begin{equation*}
\begin{split}
&\frac{\eps}{2}\int\int_{\Pi_T}\varphi(\rho_x, m_x)\nabla^2\eta^*(\rho_x, m_x)^\top dxdt\\
\leq&\int\int_{\Pi_T}C(\eps \phi_x^2m^2\rho^{-\gamma}+\|V\|_{L^{\infty}})dxdt\\
&+\int\int_{\Pi_T}\left[\eta^*\varphi_t+q^*\varphi_x+\eps\eta^*\varphi_{xx}+(\frac{m^{3}}{2\rho^{2}}
+\frac{\gamma}{\gamma-1}m\rho^{\gamma-1}p_0)\|U\|_{L^{\infty}}\varphi\right]dxdt\\
\leq& C(\varphi).
\end{split}
\end{equation*}
Thus we have arrived that
\begin{equation}\label{locestimateexclude}
 \eps\rho^{\gamma-2}\rho_x^2+\eps\rho u_x^2\in L^1_{loc}(\Pi_T).
\end{equation}
Note that when $\gamma>2,\rho^{\gamma-2}\rho^2_x$ is degenerate near $\rho=0.$ For simplicity, we assume that $1<\gamma\leq2.$
For any weak entropy-entropy flux pairs given in \eqref{2.6}, as in \eqref{4.1}, we have
\begin{equation}\label{4.3}
\begin{split}
\eta_t+q_x&=\eps\eta_{xx}-\eps(\rho_x, m_x)\nabla^2\eta(\rho_x, m_x)^\top
+(\eta_\rho U+\eta_mV)-2\eps \eta_m\rho_x\phi_x=:\sum_{i=1}^4J_i.
\end{split}
\end{equation}
By \eqref{locestimateexclude} and the boundedness of $\rho$ and $\frac{m}{\rho}$ from Theorem \ref{main3.1}, we have
\begin{equation}\label{L1}\eps\rho_x^2+\eps m_x^2\in L^1_{loc}(\Pi_T).\end{equation}
 For $J_1$,
 \begin{eqnarray*}&&\Big|\int\int_{\Pi_T} \eps\eta_{xx}\varphi dxdt\Big|\le \Big|\int\int_{\Pi_T} \eps\eta_{x}\varphi_x dxdt\Big| \le \Big|\int\int_{\Pi_T} \eps(\eta_{\rho}, \eta_m)(\rho_x, m_x)^\top \varphi_x dxdt\Big|\\
 &&\leq C\Big|\int\int_{\Pi_T} \eps(\rho_x + m_x) \varphi_x dxdt\Big|\le C\int\int_{\Pi_T}\eps^{\frac{4}{3}}(\rho_x^2 + m_x^2)+\eps^{\frac{2}{3}}\varphi_x^2 dxdt,
 \end{eqnarray*} then by (\ref{L1}), it is obvious that
$J_1$ is compact in $H^{-1}_{loc}(\Pi_T).$
Note that for any weak entropy, the Hessian matrix $\nabla^2\eta$ is controlled by  $\nabla^2\eta^*$ (see \cite{Lions2}), i.e.,
\begin{equation}\label{4.4}
(\rho_x, m_x)\nabla^2\eta(\rho_x, m_x)^\top\leq C(\rho_x, m_x)\nabla^2\eta^*(\rho_x, m_x)^\top,
\end{equation}
then $J_2$ is bounded in $L^1_{loc}(\Pi_T)$ and then compact in
$W_{loc}^{-1, \nu}(\Pi_T)$ by the embedding theorem, for some $1<\nu<2$.
Similarly, $J_3$ and $J_4$ 
are bounded in $L^1_{loc}(\Pi_T).$
Thus
$$\eta_t+q_x \text{ is compact in } W^{-1, \alpha}_{loc}(\Pi_T) \text{ for some } 1<\nu<2.$$
On the other hand,
$\eta_t+q_x \text{ is bounded in } W^{-1, \infty}_{loc}(\Pi_T).$
With the help of Lemma \ref{murat}, we conclude that
\begin{equation}\label{compact}
\eta_t+q_x~ \text{is compact in}~ H^{-1}_{loc}(\Pi_T)
\end{equation}
for all weak entropy-entropy flux pairs.
\begin{remark}
We are focusing on the uniform bound of $\rho$ and $m$. In the above argument, \eqref{compact} still holds for the case that $\gamma>2$ by a similar argument of \cite{Wang}. See also \cite{Lu}. We assume $1<\gamma\leq2$ for simplicity.
\end{remark}

\subsection{Entropy solution}
 By \eqref{compact} and the compactness framework established in
\cite{Ding, Diperna, Lions2},  we can prove that there exists
a subsequence of $(\rho^\eps,m^\eps)$ (still denoted by $(\rho^\eps,m^\eps)$) such  that
 \begin{equation}
 \label{4.6}
(\rho^\eps, m^\eps)\to(\rho, m) ~~~
  \text{ in } L^p_{loc}(\Pi_T), ~~p\geq1. \end{equation}
As in \cite{Chen2, Chen2003, Tsuge2004a, Tsuge2006},  we can prove that $(\rho, m)$ is an entropy  solution to the initial-boundary value problem \eqref{geometricc} and
$m|_{x=1}=0$
in the sense of the divergence-measure fields introduced in  \cite{Chen1999,Chen19992}. Therefore Theorem \ref{exclude} is completed.
\begin{remark}
The approach can also be applied to the Euler-Poisson system with spherical symmetry.
\end{remark}
\section{Proof of Theorem 1.2}\label{includeorigin}
\subsection{Uniform upper bound estimate}\label{upper}
Consider the system \eqref{geometric} on a cylinder $(a,b)\times \R^{+},$ with $\R^{+}=[0,+\infty),a:=a(\eps),b:=b(\eps)>1,$ and $\lim\limits_{\eps\rightarrow0}\eps a(\eps)^\delta=0$ for any $\delta\in\R$ and  $\lim\limits_{\eps\rightarrow0}b(\eps)=\infty.$ For example, $a(\eps)$ can be taken as $-\frac{1}{\ln\eps}$.
We make a scaling transformation  $$\rho=\tilde{\rho}x^c,m=\tilde{m}x^d.$$
Taking $d=(\theta+1)c>0,$ the system \eqref{geometric} can be rewritten as
\begin{eqnarray}\label{geometric222}
\left\{ \begin{split} \tilde{\rho}_t+x^{d-c}\tilde{m}_x=&-(N-1+d)x^{d-c-1}\tilde{m},\,\\ \tilde{m}_t+x^{d-c}\left(\frac{\tilde{m}^2}{\tilde{\rho}}+p(\tilde{\rho})\right)_x=&[-(2d-c+N-1)\frac{\tilde{m}^2}{\tilde{\rho}}-(2d-c)p(\tilde{\rho})]x^{d-c-1}.
\end{split}
\right.
\end{eqnarray}
If $c-d+1\neq0,$ let $\xi=\frac{1}{c-d+1}x^{c-d+1}.$ If $c-d+1=0,$ let $\xi=\ln x.$
Then \eqref{geometric222} becomes
\begin{eqnarray}\label{geometriccc}
\left\{ \begin{split} \tilde{\rho}_t+\tilde{m}_\xi=&-(N-1+d)x^{d-c-1}\tilde{m},\,\\ \tilde{m}_t+\left(\frac{\tilde{m}^2}{\tilde{\rho}}+p(\tilde{\rho})\right)_\xi=&[-(2d-c+N-1)\frac{\tilde{m}^2}{\tilde{\rho}}-(2d-c)p(\tilde{\rho})]x^{d-c-1}.
\end{split}
\right.
\end{eqnarray}
We approximate \eqref{geometriccc} by adding artificial viscosity as follows:
\begin{eqnarray}\label{geometric3}
\left\{ \begin{split} \tilde{\rho}_t+\tilde{m}_\xi=&-(N-1+d)x^{d-c-1}\tilde{m}+\eps\tilde{\rho}_{\xi\xi},\,\\ \tilde{m}_t+\left(\frac{\tilde{m}^2}{\tilde{\rho}}+p(\tilde{\rho})\right)_\xi=&[-(2d-c+N-1)\frac{\tilde{m}^2}{\tilde{\rho}}-(2d-c)p(\tilde{\rho})]x^{d-c-1}+\eps \tilde{m}_{\xi\xi}.
\end{split}
\right.
\end{eqnarray}
The initial-boundary value conditions are given as follows:
\begin{equation}\label{geometricini-vis22}
\begin{split}
(\tilde{\rho}, \tilde{m})|_{t=0}&=(\tilde{\rho}_0^\eps(x),\tilde{m}_0^\eps(x))\\
&=(\tilde{\rho}_0(x)+\eps^
{\frac{2}{\theta}}, (\frac{\tilde{m}_0(x)}{\tilde{\rho}_0(x)}+\eps)(\tilde{\rho}_0(x)+\eps^
{\frac{2}{\theta}})\chi_{[2a(\eps),b(\eps)]})\ast j^\eps, x\in [a(\eps),b(\eps)],\\
(\tilde{\rho},\tilde{m})|_{x=a(\eps)}&=(\tilde{\rho}_0^\eps (a(\eps)), \tilde{m}_0^\eps (a(\eps)))=(\tilde{\rho}_0^\eps (a(\eps)),0),\\
(\tilde{\rho},\tilde{m})|_{x=b(\eps)}&=(\tilde{\rho}_0^\eps (b(\eps)),\tilde{m}_0^\eps (b(\eps)), t>0,
\end{split}
\end{equation}
where $j^\eps$ is the standard mollifier and $\chi$ is the  characteristic function.
As in the proof of Theorem \ref{exclude}, the key point is to derive the uniform upper bound of the approximate solution $\tilde{\rho}$ and $\tilde{m}$.
By the definition of Riemann invariants, we have
$$w=\frac{m}{\rho}+\rho^\theta=\left(\frac{\tilde{m}}
{\tilde{\rho}}+\tilde{\rho}^\theta\right)x^{c\theta}:=\tilde{w}x^{c\theta}, ~z=\frac{m}{\rho}-\rho^\theta=\left(\frac{\tilde{m}}{\tilde{\rho}}-\tilde{\rho}^\theta\right)x^{c\theta}:=\tilde{z}x^{c\theta}.$$
Similarly, we have $$\lambda_1=\frac{m}{\rho}-\theta\rho^\theta=\left(\frac{\tilde{m}}{\tilde{\rho}}-\theta\tilde{\rho}^\theta\right)x^{c\theta}:=\tilde{\lambda_1}x^{c\theta},~\lambda_2=\frac{m}{\rho}+\theta\rho^\theta=\left(\frac{\tilde{m}}{\tilde{\rho}}+\theta\tilde{\rho}^\theta\right)x^{c\theta}:=\tilde{\lambda_2}x^{c\theta}.$$
It is obvious that the system \eqref{geometric3} is equivalent to the following system
\begin{eqnarray}\label{geometricvis2}
\left\{ \begin{split}
\displaystyle \rho_t+m_x=&-\frac{N-1}{x}m+\eps\left[\rho_{xx}x^{2(d-c)}+(d-3c)\rho_xx^{2(d-c)-1}\right.\\
&+\left. c(2c+1-d)\rho x^{2(d-c)-2}\right],\,\\
 m_t+\left(\frac{m^2}{\rho}+p(\rho)\right)_x=&-\frac{N-1}{x}\frac{m^{2}}{\rho}+\eps\left[m_{xx}x^{2(d-c)}-(c+d)m_xx^{2(d-c)-1}\right.\\
 &+\left.d(c+1)m x^{2(d-c)-2}\right],x\in[a(\eps),b(\eps)],t>0.
\end{split}
\right.
\end{eqnarray}

\

By the rescaled Riemann invariants $\tilde{w}$ and $\tilde{z}$, we have
\begin{eqnarray}\label{wwzz222}
\left\{ \begin{split}
\tilde{w}_t+\tilde{\lambda}_2 \tilde{w}_\xi=&\eps \tilde{w}_{\xi\xi}+2\eps \frac{\tilde{\rho}_\xi}{\tilde{\rho}}\tilde{w}_\xi-\eps\theta(\theta+1)\tilde{\rho}^{\theta-2}\tilde{\rho}_\xi^2\\
&+\left[(c-d)\frac{\tilde{m}^2}{\tilde{\rho}^2}-\theta(N-1+d)\tilde{\rho}^{\theta-1}\tilde{m}-(2d-c)\frac{\theta^2}{\gamma}\tilde{\rho}^{2\theta}\right]x^{d-c-1},\\
\tilde{z}_t+\tilde{\lambda}_1  \tilde{z}_\xi=&\eps\tilde{z}_{\xi\xi}+2\eps \frac{\tilde{\rho}_\xi}{\tilde{\rho}}\tilde{z}_\xi+\eps\theta(\theta+1)\tilde{\rho}^{\theta-2}\tilde{\rho}_\xi^2\\
&+\left[(c-d)\frac{\tilde{m}^2}{\tilde{\rho}^2}+\theta(N-1+d)\tilde{\rho}^{\theta-1}\tilde{m}-(2d-c)\frac{\theta^2}{\gamma}\tilde{\rho}^{2\theta}\right]x^{d-c-1}.
\end{split}
\right.
\end{eqnarray}
Setting the control functions $(\phi,\psi)=(M_3+2\eps, 0)$ and  by the initial and boundary data \eqref{geometricini-vis22}, we have that
 \begin{equation}\label{4.6}
\begin{aligned}
&\frac{m_0}{\rho_0}+\rho_0^{\theta} \leq (M_3 +2\eps)x^{c\theta},\,\frac{m_0}{\rho_0}-\rho_0^{\theta}\geq 0,\text{ a.e. $x\in [a(\eps),b(\eps)]$},\\
&\left(\frac{m}{\rho}+\rho^{\theta}\right)|_{x=a(\eps)}\leq (M_3+2\eps) a(\eps)^{c\theta},\left(\frac{m}{\rho}-\rho^{\theta}\right)|_{x=a(\eps)} \geq 0,\\
&\left(\frac{m}{\rho}+\rho^{\theta}\right)|_{x=b(\eps)}\leq (M_3+2\eps)  b(\eps)^{c\theta},\left(\frac{m}{\rho}-\rho^{\theta}\right)|_{x=b(\eps)} \geq 0.
\end{aligned}
\end{equation}

Define the modified Riemann invariants $(\hat{w},\hat{z})$ as
\begin{equation}\label{r1}
\hat{w}=\tilde{w}-\phi, ~~\hat{z}=\tilde{z}+\psi.
\end{equation}
The system \eqref{wwzz222} becomes
\begin{eqnarray}\label{rst}
\displaystyle\left\{ \begin{split} &\hat{w}_t+\left(\tilde{\lambda}_2-2\eps\frac{\tilde{\rho}_\xi}{\tilde{\rho}} \right)\hat{w}_\xi
=\eps\hat{w}_{\xi\xi}+a_{11}\hat{w}
+a_{12}\hat{z}+R_1,\\
&\hat{z}_t+\left(\tilde{\lambda}_1-2\eps\frac{\tilde{\rho}_\xi}{\tilde{\rho}} \right) \hat{z}_\xi
=\eps\hat{z}_{\xi\xi}+a_{21}\hat{w}
+a_{22}\hat{z}+R_2,
\end{split}
\right.
\end{eqnarray}
where
\begin{equation*}
\begin{split}
a_{11}&=0,\,a_{12}=-\theta(N-1+d)\tilde{\rho}^\theta x^{d-c-1}\leq0,\,\\
a_{21}&=0,\,a_{22}=\left[\frac{1}{4}\left[(c-d)-\theta(N-1+d)-\theta^2c\right]\hat{z}+\frac{1}{2}\left[(c-d)+\theta^2c\right]\tilde{w}\right]{x^{d-c-1}},\\
R_1&=\left[(c-d)\frac{\tilde{m}^2}{\tilde{\rho}^2}-\theta(N-1+d)\tilde{\rho}^{2\theta}-\theta^2c\tilde{\rho}^{2\theta}\right]x^{d-c-1}-\eps\theta(\theta+1)\tilde{\rho}^{\theta-2}\tilde{\rho}_{\xi}^2\leq0,\\
\,R_2&=\frac{1}{4}\left[(c-d)+\theta(N-1+d)-\theta^2c\right]\tilde{w}^2 {x^{d-c-1}}+\eps\theta(\theta+1)\tilde{\rho}^{\theta-2}\tilde{\rho}_{\xi}^2\\
&\geq\frac{1}{4} \theta(N-1)\tilde{w}^2 {x^{d-c-1}}\geq0.
\end{split}
\end{equation*}
By the
maximum principle Lemma \ref{initial maximum}, we have $$\tilde{w}(\xi,t)\leq (M_3+2\eps), \tilde{z}(\xi,t)\ge 0,$$ which implies that
$$0\le \tilde{\rho}^{\theta}(x,t)\leq \frac{M_3}{2}+\eps, 0\leq \tilde{m}(x,t)\leq (M_3+2\eps)\tilde{\rho}(x,t),$$ i.e., \begin{equation}\label{solution3}
0\leq\rho^{\theta}(x,t)\leq (\frac{M_3}{2}+\eps)x^{c\theta}, 0\leq m(x, t)\leq (M_3+2\eps)\rho(x, t)x^{c\theta}.
\end{equation}
Thus we have
 \begin{theorem}\text($L^{\infty}$ estimate: including the origin)\label{main3} Let $\gamma>1.$ Assume that for any positive constants $c$ and $M_3$, the initial and boundary data satisfy
\begin{equation}\label{ini2}
\begin{aligned}
&\rho_0(x)\ge \eps^{\frac{2}{\theta}}x^c,~\frac{m_0}{\rho_0}+\rho_0^{\theta} \leq (M_3 +2\eps)x^{c\theta},\,\frac{m_0}{\rho_0}-\rho_0^{\theta}\geq 0,\text{ a.e. $x\in [a(\eps),b(\eps)]$},\\
&\left(\frac{m}{\rho}+\rho^{\theta}\right)|_{x=a(\eps)}\leq (M_3+2\eps) a(\eps)^{c\theta},\left(\frac{m}{\rho}-\rho^{\theta}\right)|_{x=a(\eps)} \geq 0,\\
&\left(\frac{m}{\rho}+\rho^{\theta}\right)|_{x=b(\eps)}\leq (M_3+2\eps)  b(\eps)^{c\theta},\left(\frac{m}{\rho}-\rho^{\theta}\right)|_{x=b(\eps)} \geq 0,
\end{aligned}
\end{equation}
then the solution of \eqref{geometricvis2} and \eqref{4.6} satisfies
\begin{equation}\label{solution2}
0\leq\rho(x,t)\leq (\frac{M_3}{2}+\eps)^{\frac{1}{\theta}}x^c,~0\leq m(x, t)\leq (M_3+2\eps)\rho(x, t)x^{c\theta},
\end{equation}
for $x\in[a(\eps),b(\eps)]$.
\end{theorem}
\subsection{Lower bound estimate}\label{low2}
When $c-d+1\neq0$, set $\tilde{v}=\ln\tilde{\rho}$, then we get a scalar equation for $\tilde{v},$
\begin{equation}
\tilde{v}_t+\tilde{v}_\xi\tilde{u}+\tilde{u}_\xi=\eps \tilde{v}_{\xi\xi}+\eps \tilde{v}_\xi^2-\frac{N-1+b}{(c-d+1)\xi}\tilde{u}.
\end{equation}
Note that the velocity $\tilde{u}=\frac{\tilde{m}}{\tilde{\rho}}$ is uniformly bounded, i.e., $|\tilde{u}|\le C$, and  
 $$\tilde{v}(\xi,0)=\tilde{v}_0(\xi),\frac{1}{c-d+1}(a(\eps))^{c-d+1}\leq\xi\leq\frac{1}{c-d+1}(b(\eps))^{c-d+1}.$$
Following the same way as in  Subsection \ref{low}, we can show that $\tilde{\rho}\geq e^{-C(\eps, t)}.$ A similar argument can be applied to the case when $c-d+1=0.$
Consequently, we conclude the following theorem for the global existence of approximate solutions.
\begin{theorem}
Under the assumption of the previous theorem, for any time $T>0,$  there exists a positive constant   $\eps_0$  such that for $0<\eps<\eps_0$, the initial-boundary value problem
 \eqref{geometric3}-\eqref{geometricini-vis22} admits a unique classical solution on $[a(\eps), b(\eps)]\times [0, T]$ satisfying
\begin{equation*}
e^{-C(\eps, T)}x^c\leq \rho(x, t)\leq (\frac{M_3}{2}+\eps)^{\frac{1}{\theta}}x^c, ~~0\leq m(x, t)\leq (M_3+2\eps)\rho(x, t)x^{c\theta}.
\end{equation*}
\end{theorem}

\subsection{$H^{-1}_{loc}$ compactness of entropy pair}\label{entropy}
For any $T\in(0, \infty)$, let $\Pi_T=(0,+\infty)\times (0, T)$.
Let $K\subset\Pi_T$ be any compact set, and choose $\varphi\in C_c^\infty(\Pi_T)$ such that $\varphi|_{K}=1,$ and $0\leq\varphi\leq1$. When $\eps$ is small, $K\subset (a(\eps),b(\eps))\times (0, T)$.
Similarly as in subsection \ref{entropy11}, multiplying \eqref{geometricvis2} by  $\nabla\eta^* \varphi$ with $\eta^*$ the mechanical entropy,
we obtain
\[
(\rho_x, m_x)\nabla^2\eta^*(\rho_x, m_x)^\top x^{2(d-c)}=(p_0\gamma\rho^{\gamma-2}\rho_x^2
+\rho u_x^2)x^{2(d-c)},
\]
and
\begin{equation}\label{locestimate}
 (\eps\rho^{\gamma-2}\rho_x^2+\eps\rho u_x^2)x^{2(d-c)}\in L^1_{loc}(\Pi_T),
\end{equation}
which implies that for any weak entropy-entropy flux pairs $(\eta,q)$ given in \eqref{2.6}, it holds
\begin{equation}\label{4.15}
\eta_t+q_x~ \text{is compact in}~ H^{-1}_{loc}(\Pi_T).
\end{equation}
Since the proof is almost the same as in Subsection \ref{entropy11}, we omit it here.
\subsection{Entropy solution}\label{strong}

 By \eqref{4.15} and the compactness framework established in
 \cite{Ding, Diperna, Lions2},  we can prove that there exists
 a subsequence of $(\rho^\eps,m^\eps)$ (still denoted by $(\rho^\eps,m^\eps)$) such  that
 \begin{equation}
 \label{4.16}
 (\rho^\eps, m^\eps)\to(\rho, m) ~~~
 \text{ in } L^p_{loc}(\Pi_T), ~~p\geq1. \end{equation}

Note that $\rho=O(1)x^c$ and $m=O(1)x^{d}$ so that the  right hand sides of \eqref{geometric}, that is, $\frac{1}{x}m$ and $\frac{1}{x}\frac{m^2}{\rho}$ are integrable near the origin with respect to $x$.  As in \cite{Chen2, Chen2003, Tsuge2004a, Tsuge2006},  we can prove that $(\rho, m)$ is an entropy  solution to the problem \eqref{geometric} and the test function $\Phi(x,t)$ can contain the origin.  Therefore, the proof of Theorem \ref{include} is completed. \\

It is remarked that the entropy solution obtained above is exactly the entropy solution to the Cauchy problem of isentropic gas dynamics system with spherical symmetry.

\bigskip

\section{Appendix}
Proof of Lemma \ref{density}: 
\begin{proof}
 We decompose $w$ into $\displaystyle w=\sum_{i=0}^3w_i$, where $w_2$ and $w_3$ are
\begin{equation}
\begin{split}w_2(x,t)=\int_0^t\int_a^b\Gamma(x-\xi,t-\tau)f_2(\xi,\tau)d\xi d\tau,\\ w_3(x,t)=\int_0^t\int_a^b\Gamma_x(x-\xi,t-\tau)h(\xi,\tau)d\xi d\tau,\end{split}\end{equation}
and $\Gamma$ is the heat kernel:
\begin{equation}\label{gamma}
\Gamma(x-\xi,t-\tau)= \begin{cases}
\frac{1}{(4\pi\eps(t-\tau))^{\frac{1}{2}}}e^{-\frac{|x-\xi|^2}{4\eps(t-\tau)}},&t>\tau;\\
0,&t\leq \tau,
\end{cases}
\end{equation}
and $w_0$ and $w_1$ are the solutions of the following problems respectively:
\begin{eqnarray}\label{0}
(P_0):\left\{ \begin{array}{ll}
\displaystyle w_0-\eps w_{0xx}=0, a<x<b, \,\\
\displaystyle w_0|_{t=0}=\varphi(x), a<x<b,\\
\displaystyle w_0|_{x=a}=\varphi(a)-w_2(a,t)-w_3(a,t),w_0|_{x=b}=\varphi(b)-w_2(b,t)-w_3(b,t),
\end{array}
\right.
\end{eqnarray}
\begin{eqnarray}\label{1}
(P_1):\left\{ \begin{array}{ll}
\displaystyle w_{1t}-\eps w_{1xx}=f_1(x,t), a<x<b,\,\\
\displaystyle w_1|_{t=0}=0, a<x<b,\\
\displaystyle w_1|_{x=a}=0,w_1|_{x=b}=0.
\end{array}
\right.
\end{eqnarray}
It is obvious that $w_{2t}-\eps w_{2xx}=f_2(x,t)$ and $w_{3t}-\eps w_{3xx}=h_x(x,t)$. Thus $w= \sum_{i=0}^3w_i$ is the unique solution of the equation \eqref{u}.
Note that
\begin{equation}\label{gamma2}
\begin{split}
0\leq\Gamma(x-\xi,t-\tau)\leq\frac{C}{(t-\tau)^{\frac{1}{2}}},\\
|\Gamma_x(x-\xi,t-\tau)|\leq\frac{C}{(t-\tau)^{\frac{3}{2}-\alpha}|x-\xi|^{2\alpha-1}},
\end{split}
\end{equation}
for any $\frac{1}{2}<\alpha<1$. For example  $\alpha=\frac{2}{3}$, then $ \Gamma_x$ is integrable with respect to $\xi$ and $\tau$ up to time $T>0$.
From \eqref{gamma2}, we obtain that $w_2$ and $w_3$ are bounded. Moreover, they are also bounded on the boundary $x=a,b$.
We then turn to the problem $(P_0) $ and find that its solution $w_0$  is bounded. For the problem $(P_1)$, the solution $w_1\geq0$ due to $f_1(x,t)\geq0.$
Therefore, the proof of Lemma \ref{density} is completed.
\end{proof}

\section*{Acknowledgments}
Feimin Huang was partially supported by National Center for Mathematics and Inter- disciplinary Sciences, AMSS, CAS and NSFC Grant No. 11371349 and 11688101. Tianhong Li was
supported by National Natural Science Foundation of China  under contracts 10931007 and 11771429. Difan Yuan was supported by China Scholarship Council No.201704910503. The authors would like to thank the anonymous referees for valuable comments and suggestions.
%

\end{document}